%% file: main.tex
\documentclass[conference]{IEEEtran}
\IEEEoverridecommandlockouts
\ifCLASSINFOpdf
\else
\fi
\usepackage[cmex10]{amsmath}
\usepackage{graphicx}
\usepackage{acronym}
\usepackage{xspace}

\begin{document}

%
\title{Preliminary Results in Current Profile Estimation and Doppler-aided Navigation\\for Autonomous Underwater Gliders}
%
\author{
Jonathan Jonker$^{1}$,
Andrey Shcherbina$^{2}$,
Richard Krishfield$^{3}$,
Lora Van Uffelen$^{4}$,\\
Aleksandr Aravkin$^{5}$,
Sarah E. Webster$^{2*}$
\thanks{$^{1}$Department of Mathematics, University of Washington, Seattle, WA USA}
\thanks{$^{2}$Applied Physics Laboratory, University of Washington, Seattle, WA USA}
\thanks{$^{3}$Woods Hole Oceanographic Institution, Woods Hole, MA USA}
\thanks{$^{4}$Ocean Engineering, University of Rhode Island, Narragansett, RI USA}
\thanks{$^{5}$Applied Mathematics, University of Washington, Seattle, WA USA}
\thanks{*Corresponding author: {\tt\small swebster@apl.washington.edu}}%
}

\maketitle

\vspace{-0.2in}

\input{acronym_list}
\input{commands}

\begin{abstract}
    
This paper describes the development and experimental results of navigation algorithms for an autonomous underwater glider (AUG) that uses an on-board acoustic Doppler current profiler (ADCP). AUGs are buoyancy-driven autonomous underwater vehicles that use small hydrofoils to make forward progress while profiling vertically.
During each dive, which can last up to 6 hours, the Seaglider AUG used in this experiment typically reaches the depth of 1000 m and travels 3-6 km horizontally through the water, relying solely on dead-reckoning. Horizontal through-the-water (TTW) progress of AUG is 20-30 cm/s, which is comparable to the speed of the stronger ocean currents. Underwater navigation of an AUG in the presence of unknown advection therefore presents a considerable challenge. We develop two related formulations for post-processing. Both use ADCP observations, through the water velocity estimates, and GPS fixes to estimate current profiles. However, while the first solves an explicit inverse problem for the current profiles only, the second solves a deconvolution problem that infers both states and current profiles using a state-space model.

Both approaches agree on their estimates of the ocean current profile through which the AUG was flown using measurements of current relative to the AUG from the ADCP, and estimates of the AUGs TTW velocity from a hydrodynamic model.  The result is a complete current profile along the AUGs trajectory, as well as over-the-ground (OTG) velocities for the AUG that can be used for more accurate subsea positioning.
Results are demonstrated using 1 MHz ADCP data collected from a Seaglider AUG deployed for 49 days off the north coast of Alaska during August and September 2017. The results are compared to ground truth data from the top 40 meters of the water column, from a moored, upward-facing 600 kHz ADCP.
Consequences of the state-space formulation are discussed 
in the Conclusions section. 
\end{abstract}

\section{Introduction}

\begin{figure}[ht]
  \centering
  \includegraphics[width=\columnwidth]{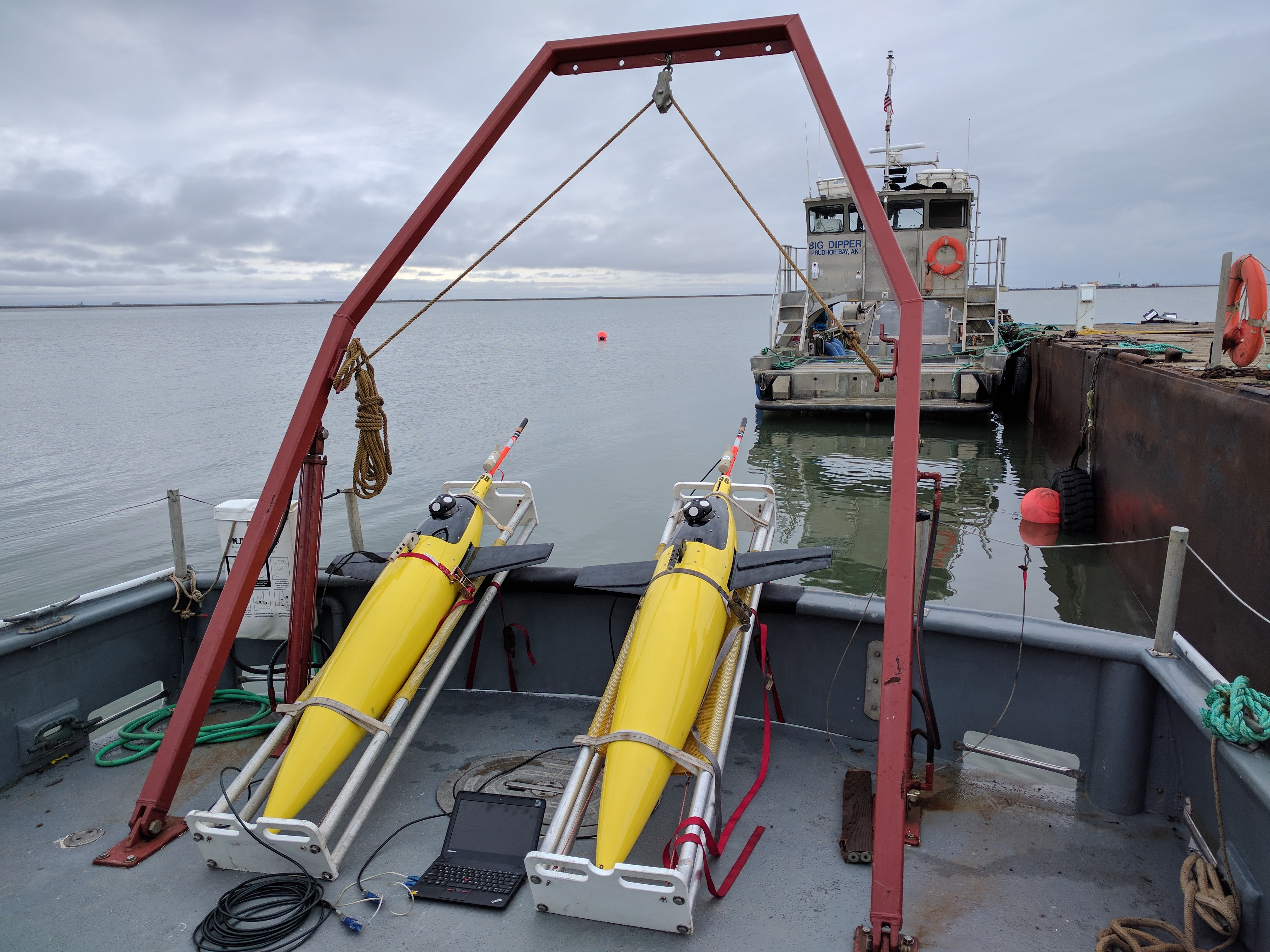}
  \caption{Ready for launch, Seagliders SG196 and SG198 are loaded on the R/V Ukpik in Prudhoe Bay, AK, with the upward-facing ADCPs visible where they are installed in the aft fairing.}
  \label{fig:SG}
  \vspace{-0.2in}
\end{figure}

The goal of this work is to simultaneously estimate a) absolute (Earth-referenced) ocean velocity profile, and b) the absolute autonomous underwater glider (AUG) path over the bottom, based on on-board acoustic Doppler current profiler (ADCP) observations of relative current velocity profiles. Figure \ref{fig:SG} shows two Seagliders with upward-facing ADCPs mounted in the aft fairing.

The main challenge of AUG-based ADCP velocity profiling is that each velocity profile is measured relative to the through-the-water (TTW) motion of the glider, as opposed to the earth referenced glider velocity. The AUG's TTW velocity can be inferred from the AUG dynamic model, but is subject to uncertainty because the model is based on steady state flight and doesn't take roll into account. This lack of a georeferenced platform velocity requires a simultaneous estimation of current and glider TTW velocity using all available data---ADCP velocity profiles, surrounding water density, glider buoyancy engine state, glider attitude and angle of attack, glider depth, and GPS positions at the start and end of the dive.

Here we describe two different frameworks for performing this estimation---a linear global inverse, and a flexible state space approach. 
In this initial work, we use a pre-computed estimate of the glider TTW velocity from the hydrodynamic model for both methods. The overarching future goal is to use the flexibility 
of the state-space model to include the nonlinear hydrodynamic model and nonlinear range measurements as part of the approach.

\section{Background}

Historically, ADCP measurements were made from ships with with relatively low frequency systems (70 kHz) or from moorings with low or mid-frequency instruments (300 kHz). The 300 kHz ADCPs are also common on large underwater vehicles, in particular instruments that can be used as a Doppler velocity log to provide ``bottom lock''---i.e., vehicle velocity relative to the seafloor \cite{yoerger-1998-ABE} or, in some specialized cases, the underside of ice \cite{mcfarland-2015-icerelnav}.

To provide better resolution of deep currents from shipboard measurements, lowered ADCP methods were developed using overlapping shear traces from a higher frequency (typically 150 or 300 kHz) ADCP lowered from a (relatively) stationary ship.  The shear traces are then used to reconstruct the full current profile. The two main lowered ADCP processing frameworks are the shear method, originally described by \cite {firing-1990-ladcp,fischer-1993-ladcp}, and the velocity inversion method \cite{visbeck-2002-ladcp}.
%

In the last 10 years, there has been increased interest in using ADCPs from autonomous underwater vehicles to characterize the currents between the surface and the seafloor and to close the gap of uncertainty in vehicle drift between when GPS is available at the surface and bottom lock is obtained within range of the seafloor \cite{mstanway-2010a,stanway-2011-adcpnav,medagoda10_oceans}.

Aided by the development of very high frequency ADCPs (1 MHz) specifically for autonomous underwater gliders, the two main lowered ADCP processing frameworks---the shear method and the inversion method---have been modified to estimate depth-varying current profiles from gliders by \cite{thurnherr-2015-slocum-adcp} and \cite{todd-2011-JGRC,todd-2017-JAOT} respectively. A related, but distinctly different method based on non-linear objective mapping of ADCP shear has been also developed by \cite{davis-2010-adcpglider}.
%
%
In all implementations (\cite{thurnherr-2015-slocum-adcp,todd-2017-JAOT,davis-2010-adcpglider}), the authors needed to contend with the limitation of Slocum and Spray glider ADCPs that only collect data during one half of the dive (either ascending or descending).  

%
Our contribution builds on and extends the velocity inversion method with several key innovations. We use current data from the entire dive (as opposed to just during descent), which is crucial for multi-hour dives where, as we show, the current profiles on descent versus ascent differ significantly. We then explicitly separate the over-the-ground glider velocity into the drift velocity (as a result of advection), and the glider's horizontal through-the-water velocity. This enables us to directly incorporate hydrodynamic model velocity estimates and independently control the smoothness regularization of the different velocity components.

\section{Finding the Current Profile by Inversion}
\label{sec:first}

In this section we develop an approach to use 
(1) on-board ADCP observations of relative current velocity profiles, (2) dive start and end GPS coordinates and (3) a measure of through the water (TTW) velocity obtained from the glider to simultaneously estimate (a) absolute (Earth-referenced) ocean velocity profile, and (b) the glider path through the water.
We present the derivation for a scalar velocity $v$, which can be understood to be either one of the velocity components $(u,v)$, or a complex-number representation of velocity, $u+iv$.

\subsection {Formulation of the inverse problem}
Similar to \cite{visbeck-2002-ladcp,todd-2011-JGRC}, glider-borne ADCP observations of horizontal currents are considered to be a sum of three unknown components: 
$$u^a(z,t)=u^o(z)-u^g(t)+\epsilon(z,t),$$
where 
\begin{itemize}
\item $u^o$ is the absolute ocean velocity 
\item $u^g$ is the over-the-ground (OTG) velocity of the glider platform 
\item $\epsilon$ is the ADCP measurement noise. 
\end{itemize}
Here, we separate the OTG glider velocity into the drift velocity $u^d=u^o(z^g(t))$, equal to the ocean velocity at the glider depth $z^g$, and the glider horizontal TTW "propulsion" speed $u^p$, which is determined by the glider control and its flight dynamics:
\[
u^g(t) = u^o(z^g(t)) + u^p(t).
\]
Thus the inverse problem is to obtain $u^o(z)$ and $u^g(t)$ such that the discrepancy with the observed values of $u^a$ are minimized,
\begin{equation}
\label{eq:min}
\min{\|u^a(z,t)-(u^o(z)-u^o(z^g(t))-u^p(t))\|_2^2}.
\end{equation}
This formulation uses the least squares loss, corresponding
to a Gaussian assumption on the noise $\epsilon$. 
More robust losses are of interest and will be studied in future work. 

This formulation is more complex than that of \cite{todd-2011-JGRC}, as it requires interpolation of ocean velocity onto the glider location ($u^o(z^g(t))$), 
but it allows explicit treatment of the glider state. Additionally, applying smoothness regularization  to $u^p$ 
as discussed below is more appropriate than to $u^g$. 

\subsection{Discrete formulation}
We consider a set of individual ADCP observations $u^a_{ij}=u^a(z_{ij},t_j)$ at depth cells $z_{ij}$ and times $t_j$, with the
 $i\in [1,I]$ being the ADCP range cell index, and $j\in [1,J]$ being the sample time (ping) index.  All the valid observations form an observation column vector 
$$u^a=[u^a_1...u^a_k]^\top,$$
where the index $k\in [1, K], K\le IJ$ enumerates valid observations. The vectors $[z_k]$ and $[t_k]$ represent the corresponding depth and time coordinates of the $k$-th valid observation.

The unknown state vector 
$$x=\begin{bmatrix}
{u^g} \\ u^o
\end{bmatrix}$$
consists of a vector of the unknown glider OTG velocities 
$${u^g}= [u^g_1 ... u^g_M]^\top$$ 
defined on a temporal grid $\hat{t}_m, m\in[1,M]$, and a vector of the unknown ocean velocities 
$$u^o=[u^o_1 ... u^o_L]^\top$$
defined on a regular vertical grid $\hat{z}_l=(l-1)\Delta \hat{z}, l\in [1,L]$. For convenience, the temporal grid $\hat{t}_m$ is taken to be a superset of the sample times $t_j$, with the extra intervals added during the gaps in the ADCP record (at the beginning and the end of the dive, as well as during a brief intermission at the bottom of the dive); the spacing of the extra intervals is set to the average ADCP sampling interval.
Discrete formulation of the ADCP sampling problem \eqref{eq:min} then becomes 
$$
u^a= - H^t u^g + H^z u^0 + \epsilon, 
$$
where $\epsilon$ is the noise vector.
The matrix operator $H^t$ represents temporal interpolation from the time grid $\{\hat{t}_m\}$ onto the sampling times $\{t_k\}$, and 
$H^{z}$ represents spatial interpolation from the vertical grid $\{\hat{z}_l\}$ onto $\{z_k\}$.
Since $\{t_k\}\subset \{\hat{t}_m\}$ by construction, $H^t$ is simply a subsampling matrix,
$$
H^t_{km}=\begin{cases} 1,&\text{if }t_k = \hat{t}_m\\0,&\text{otherwise } \end{cases}.
$$
The formulation of vertical interpolation matrices depends on the chosen interpolation method. Visbeck \cite{visbeck-2002-ladcp} method is equivalent to the nearest-neighbor interpolation. Here, we employ linear interpolation, corresponding to a matrix operator
$$
H^z_{kl}=\begin{cases}
1-(\hat{z}_l-z_k)/\Delta \hat{z},&\hat{z}_{l-1} \le z_k <\hat{z}_l\\
1-(z_k-\hat{z}_l)/\Delta \hat{z},&\hat{z}_l \le z_k <\hat{z}_{l+1}\\
0,& \text{otherwise.}
\end{cases}
$$

\subsection {Additional Regularization}
\label{sec:inverse.constraint}

The inversion framework balances physical knowledge of the glider's motion (such
as start- and end-of-dive GPS fixes) against assumptions about the physical environment, including smoothness of current profiles in time and space.

\textbf{Start- and End-of-Dive GPS:} GPS fixes before and after the dive provide tie-points on the integral $\int u^gdt$ over the duration of the dive. The discrete formulation of these relationships is given by 
$$
Sx\approx s,
$$
where $s$ is the scalar (complex) horizontal displacement between the two GPS fixes,
\(
S=\begin{bmatrix}w & 0\end{bmatrix},
\)
and $w$ represents the weights corresponding to trapezoidal rule integration, 
$$ w_m=\begin{cases}
0.5(\hat{t}_2-\hat{t}_1), &m=1\\
0.5(\hat{t}_{m+1}-\hat{t}_{m-1}), &1<m<M\\
0.5(\hat{t}_M-\hat{t}_{M-1}), &m=M
\end{cases}
$$

\textbf{Glider Velocity:} A measure of the glider's TTW velocity, such as that provided by the glider hydrodynamic model, can be used to provide a constraint on $u^p$(t) in \eqref{eq:min}. Discrete formulation of this constraint is 
\[
u^p = -H^t u^g + H^{z0} u^o + \epsilon, 
\]
where $u^{p}$ is the TTW velocity  predicted by the hydrodynamic model, and $H^{z0}$ represents spatial interpolation of ocean velocity profile onto the glider depth 
(so that  $H^{z0}  x$ is the glider drift speed); this matrix is constructed in the same way as $H^{z}$. The relative confidence, within the model, of these measurements compared to the ADCP measurements can be controlled by adjusting the relative magnitude of the noise vectors $\epsilon$ for the two measurements. As is demonstrated and discussed in the Results section, however, understanding the subtleties of how to weight the constraints is still a work in progress.

\textbf{Smoothness Regularization:}
\label{sec:inv.reg}
Additionally, we require both the ocean and the glider velocities to be smooth, which is equivalent to minimization of ${D_2u^g}$, and ${D_2u^o}$ where $D_2$ are the second derivative operators of the appropriate sizes,
$$D_2=\begin{bmatrix}
1 &-2 &1 &0 &\cdots &0\\
0 &1 &-2 &1 &\\
\vdots & & & \ddots &\ddots\\
0 &\cdots &0 &1 &-2 &1
\end{bmatrix}.
$$

\textbf{Subtleties of Regularization:}
Ideally, we want to solve an inverse problem that is aware of both the
vehicle dynamics and control and the ADCP observations, because this is the
only way to make a self-consistent estimate. The challenge, as mentioned above
with the glider velocity estimates in particular, is that measurements are inherently
noisy, often biased, and optimally choosing the relative weighting between
measurements is still a work in progress.
%
%

In the extreme case, if we had a perfect hydrodynamic model, we would not need
the inverse at all, as we would know the TTW velocity at any moment (based
solely on the buoyancy and pitch control), and therefore would have absolute
ocean velocity estimates from every trace. At the other extreme, if the ADCP
measurements were perfect and extended close to the vehicle, we would have a
direct measure of the TTW velocity from that information alone and would not
need the dynamic model.

In reality, the ADCP is noisy and has gaps (notably, the blanking
distance). And without the hydrodynamic model, we can't distinguish between the legitimate accelerations due to AUG control (which
should be kept) and spurious accelerations arising from ADCP noise (which should
be eliminated)---and we would be forced to either keep both, or eliminate both,
which is suboptimal.





\subsection {Least-Squares Formulation}
The final least-squares problem formulation is given by 
$$
\min_x{\|{Gx-d}\|_2^2},
$$
where
$$
G=
\begin{bmatrix}
-H^t  &H^z \\
w & 0\\
 -H^t & H^{z0}\\
r_o D_2 & 0 \\
0  & r_gD_2
\end{bmatrix}, \qquad 
d=\begin{bmatrix}{u^a}\\ s\\  u^p \\ 0 \\ 0  \end{bmatrix},
$$
and $r_o$ and $r_g$ are regularization parameters.
The solution is given by
$$
x=(G^\top G)^{-1}G^\top d,
$$
and can be computed in a more efficient manner (e.g., using a QR decomposition aware of the block structure of $G$).
\subsection {Modification: two-profile solution}
The ocean velocity profile is expected to change over the duration of the dive. Therefore, it may be reasonable to seek \emph{two} ocean velocity profiles $u^d$ and $u^u$, corresponding to the descent and asent. The state vector and the observation matrix then become
$$x=\begin{bmatrix}
{u^g} \\ 
u^d \\
u^u 
\end{bmatrix},$$
$$
H= \begin{bmatrix}
&-H^t&|&H^{zd} &|&H^{zu} &\end{bmatrix},
$$
where the interpolation matrices $H^{zd}$ and $H^{zu}$ are constructed as $H^{z}$ before, except that only those rows of $H^{zd}$ that correspond to the downcast $z_k$ are non-zero, and vise versa. 

An additional constraint is necessary, requiring the ocean velocity profiles to match at the bottom end, i.e.
$$
\begin{bmatrix}
0 & e_n^\top & -e_n^\top
\end{bmatrix}x=0.
$$
where $e_n$ is the elementary vector with $1$ in the last entry. 

The modified system is given by 

$$
\min_x{\|{Gx-d}\|_2^2},
$$
where
$$
G=
\begin{bmatrix}
-H^t  &H^{zd} & H^{zu} \\
w & 0 & 0\\
 -H^t & H^{z0d} &H^{z0u} \\
r_o D_2 & 0 & 0\\
0  & r_gD_2 & r_gD_2 
\end{bmatrix}, \qquad 
d=\begin{bmatrix}{u^a}\\ s\\  u^p \\ 0 \\ 0  \end{bmatrix},
$$

For shorter dives, assuming identical current profiles on ascent and descent may be a reasonable assumptions. For the data collected over multi-hour dives, however, as we show, there can be significant deviation in the current profile during descent and ascent, and we use the two-profile solution.

\section{Deconvolving Glider State and Current Profile: 
State-Space Approach}
\label{sec:nonlin}

In this section, we propose a state-space model that fuses available information (ADCP current profiles, GPS coordinates, and measure of smoother velocities) to simultaneously infer 
the states of the glider together with current profile estimation. The formulation uses identical information to that used in Section~\ref{sec:first}. 
However, the main innovation here is to use a more detailed navigation model, with a view toward simultaneously solving the  current mapping and glider localization problems. A direct consequence of this is an updated estimate of the state variables of the glider, along with the estimates of the current profile.

\subsection{Glider State-Space Model}
We start by constructing a linear four-dimensional state-space model for the glider positions and their derivatives (OTG velocity): 
\begin{equation}
    \label{eq:KSmodel}
\begin{aligned}
x_k &:= 
\begin{bmatrix}e^g_k & n^g_k &  \dot{e}^g_k & \dot{n}^g_k \end{bmatrix}^T, \\
x_{k+1} & = G_kx_k + \epsilon_k, \quad \epsilon_k \sim\mathcal{N}(0,Q_k) \\
z_k &= H_kx_k + \nu_k, \quad \nu_k \sim \mathcal{N}(0,R_k), \quad \text{for } k=1,N\\
 \quad G_k &= \begin{bmatrix}1 & 0 & \Delta t_k & 0\\ 0 & 1 & 0 & \Delta t_k &\\ 0 & 0 & 1 & 0\\ 0 & 0 & 0 & 1\end{bmatrix}, 
H_k = \begin{bmatrix}1 & 0 & 0 & 0 \\ 
0 & 1 & 0 & 0\end{bmatrix}
\end{aligned}
\end{equation}
The four states are east/north over-the-ground (OTG) velocities and their approximate integrals. The only measurements used here are the GPS position fixes at the beginning and end of the dive. A separate term is added later to compare OTG velocities with (measured) TTW velocities.

The state-space model~\eqref{eq:KSmodel} does not force smoothness of the velocities, unlike the previous section. 
The full-state block bi-diagonal process-discrepancy matrix 
$G$, and block diagonal matrices $H, Q, R$ are  
\[
G = \begin{bmatrix}I & 0 & & \\ -G_2 & I & \ddots & \\  &\ddots & \ddots & 0\\  &  & -G_N & I\end{bmatrix}, \quad H =\begin{bmatrix} H_1 & 0 & &\\ 0 & \ddots & \ddots & \\ & \ddots & \ddots & 0\\ & & 0 & H_N\end{bmatrix}
\]
\[
Q = \begin{bmatrix} Q_1 & 0 & &\\ 0 & \ddots & \ddots & \\ & \ddots & \ddots & 0\\ & & 0 & Q_N\end{bmatrix}, \quad R = \begin{bmatrix} R_1 & 0 & &\\ 0 & \ddots & \ddots & \\ & \ddots & \ddots & 0\\ & & 0 & R_N\end{bmatrix}
\]
We also introduce notation for the full state $X$, measurements $z$, and initial state $w$: 
\[
 x=\begin{bmatrix} x_1 \\ x_2 \\ \vdots \\ x_N \end{bmatrix}, \quad z = \begin{bmatrix} z_1 \\ z_2 \\ \vdots \\ z_N \end{bmatrix}, \quad w = \begin{bmatrix}x_0 \\ 0 \\ \vdots \\ 0\end{bmatrix}
 \]
The Kalman smoother estimate for the velocities (and their derivatives) from direct observations would be obtained by solving the single least squares problem 
\[
\min_{x}||Gx-w||^2_{Q^{-1}} + ||Hx-z||^2_{R^{-1}}
\]
\subsection{ADCP Measurements}

We now consider unknown current profiles $c_u$ and $c_v$, 
which connect to the glider states in the previous sections 
via the simple equations 
\[
c_u = \dot{e}^g + z^o_u + \epsilon_{cu}
\]
\[
c_v = \dot{n}^g + z^o_v + \epsilon_{cv}
\]
The additional decision variables $c_u,c_v$ are indexed by depth. The model linking the ADCP observations to current profiles and states is given by 
\[
\begin{aligned}
z^o_u &= A_cc_u - A_ux + \epsilon_{cu} \\
z^o_v & = A_c c_v - A_v x + \epsilon_{cv}
\end{aligned}
\]
where $A_u,A_v$ select out the OTG velocities at the time the measurements were taken and $A_c$ selects out the depths for the given measurement. 

Depending on the discretization of $c_u,c_v$, there may be  depths without associated measurements. 
To estimate currents at these depths and to smooth out the final ocean velocity estimate we impose regularization terms on $c_u$ and $c_v$:
\[
R(c_u, c_v) := \sum_{i=1}^{N-1} ({c_u}_i-{c_u}_{i+1})^2 + ({c_v}_i - {c_v}_{i+1})^2,
\]
which can be written as
\[
R(c_u, c_v) = ||A_rc_u||^2 + ||A_rc_v||^2
\]
where $A_r$ compute adjacent differences.

\subsection{Comparing OTG and TTW Velocities}

The difference between OTG and TTW Velocities depends on the current,
so we add a least squares term for this estimation. 
The three variables are connected via the equations
\[
\dot{e}^g = e\_ttw^g + c_u + \epsilon_ev
\]
\[
\dot{n}^g = n\_ttw^g + c_v + \epsilon_nv
\]
These can be written using existing variables as follows
\[
A_{OTGe}x - (A_{TTWe}x_{TTW} + A_{cvelu}c_u) + \epsilon_ev = 0
\]
\[
A_{OTGn}x - (A_{TTWn}x_{TTW} + A_{cvelv}c_v) + \epsilon_nv = 0,
\]
where $A_{OTG},A_{TTW}$ select the appropriate velocities at each time and $A_{cvel}$ selects  the current for the depth at the given time.

\subsection{Joint Inversion}

Combining the state-space model with the ADCP observations, current profiles, and velocity comparisons gives the least squares problem 
\begin{equation}
\label{eq:KSjoint}
\begin{aligned}
\min_{x,x_{TTW},c_u,c_v} 
& ||Gx-w||^2_{Q^{-1}} +||Hx-z||^2_{R^{-1}}\\
+ & \eta_1 ||A_cc_v - A_vx - z^o_v||^2 + \eta_1 ||A_cc_u - A_ux - z^o_u||^2\\
+ & \eta_2 ||A_{OTGe}x - (A_{TTWe}x_{TTW} + A_{cvelu}c_u)||^2\\
+ & \eta_2 ||A_{OTGn}x - (A_{TTWn}x_{TTW} + A_{cvelv}c_v||^2\\
+ & \eta_3 ||A_rc_u||^2 + \eta_3 ||A_rc_v||^2
\end{aligned}
\end{equation}


The joint inverse problem is interesting because it violates
the classic Kalman smoothign block tridiagonal structure. 
In particular, while $G^TG$ is sparse and block tridiagonal,  
and $H^TH$ is sparse and block diagonal, the matrix $A_v^TA_v$ is a generic 
sparse matrix. 
In the experiments, we exploit the sparsity 
of the final least squares problem~\eqref{eq:KSjoint} to 
solve the problem efficiently. 
We leave further structure-exploiting innovations for ADCP-informed navigation to future work.

\section{Experimental Data Collection}
 
\begin{figure}
  \centering
  \includegraphics[width=0.9\columnwidth]{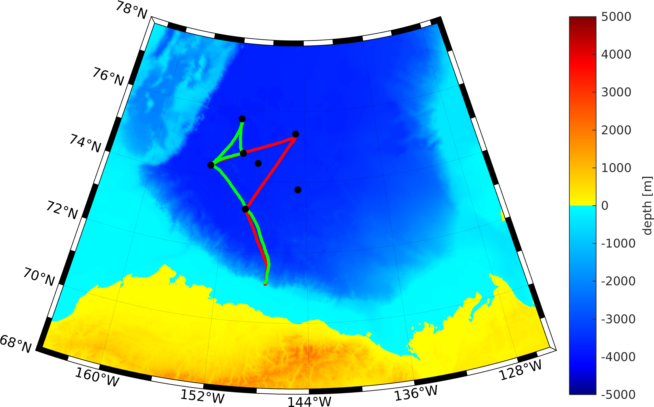}
  \caption{SG196 (green) and SG198 (red) were deployed at the shelf break north of Prudhoe Bay, AK.  From there, they flew up to and around the CANAPE mooring array (black dots) for 49 days until they were recovered by the USCGC Healy.}
  \label{fig:track}
\end{figure}

 
 As part of the Canada Basin Glider Experiment (CABAGE), two Seagliders, SG196 and SG198, were deployed on 6 August 2017 at the shelf break north of Prudhoe Bay, AK. From there, they flew up to and around the CANAPE mooring array until they were recovered on 17 September 2017, for a total of 49 days.  Together the gliders covered approximately 1730 km over the course of 712 dives, with SG196 diving to 480 m depth and SG198 diving to 750 m.  Figure \ref{fig:track} shows the glider tracklines for both a short test deployment in 2016 and the 2017 deployment.

Each glider was equipped with a Nortek Signature1000 1 MHz ADCP, as well as the standard suite of conductivity temperature (CT) sensor, pressure sensor, WHOI MicroModem, and custom-built passive marine acoustic recorders (PMARs). Figure \ref{fig:SG} shows the gliders loaded on the R/V Ukpik, ready for launch, with the upward-facing ADCPs, installed in the tail section, clearly visible.


For clarity during the discussion here, when referring to the ADCP data, we will use \emph{trace} to refer to an individual profile collected by the ADCP, and \emph{profile} to refer to the result of the inverse, i.e. the current profile for the entire dive. The ADCPs were programmed to collect a trace every 15 seconds with 2.0 m bins. Each trace typically covered 10-15 m depth, with the actual usable range varying with the amount of acoustic scatterers in the water. In practice any given depth bin was covered by 5-7 different traces.  Figure \ref{fig.snippet} illustrates a section of the current profile with overlapping ADCP traces after alignment.
 

\section{Results}

In this section we compare the results obtained using the approaches detailed in Sections~\ref{sec:first} and~\ref{sec:nonlin}. 

\begin{figure}[t]
  \includegraphics[width=\columnwidth]{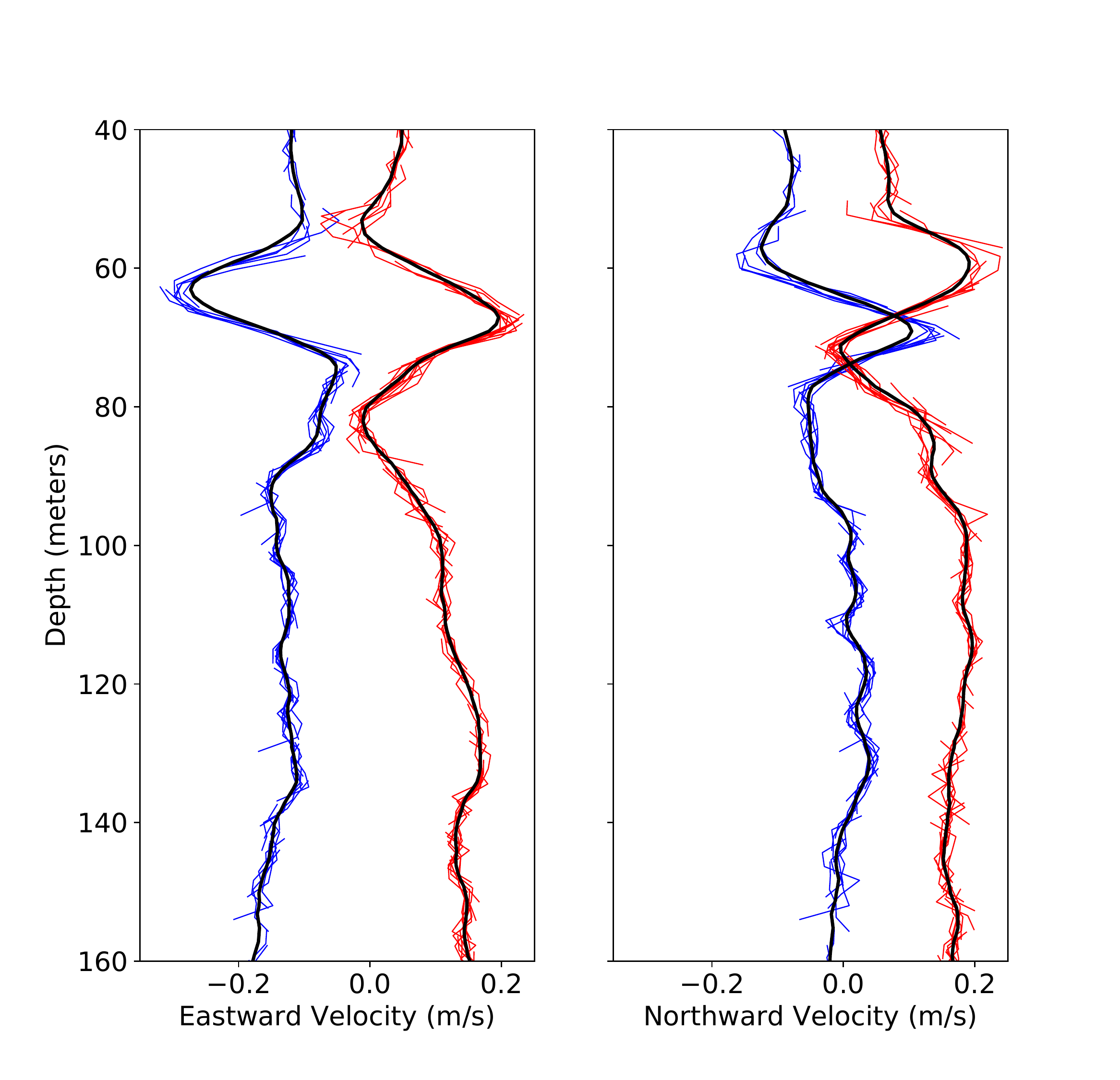}
  \caption{Overlapping ADCP traces during descent (blue) and ascent (red) for
    dive 99 of sg198, after alignment. The current profile produced by the state-space approach is shown as the thick black line.}
  \label{fig.snippet}
\end{figure}

Figure \ref{fig.snippet} shows a short section of the current profile with the individual overlapping traces that have been aligned and averaged to produce the final current profile (black).  This is from a 750 m dive, so the sections of data shown were collected 3.5 hrs apart. The difference in the current profile during ascent and descent are clearly visible.

\begin{figure}
  \centering
  \includegraphics[width=\columnwidth]{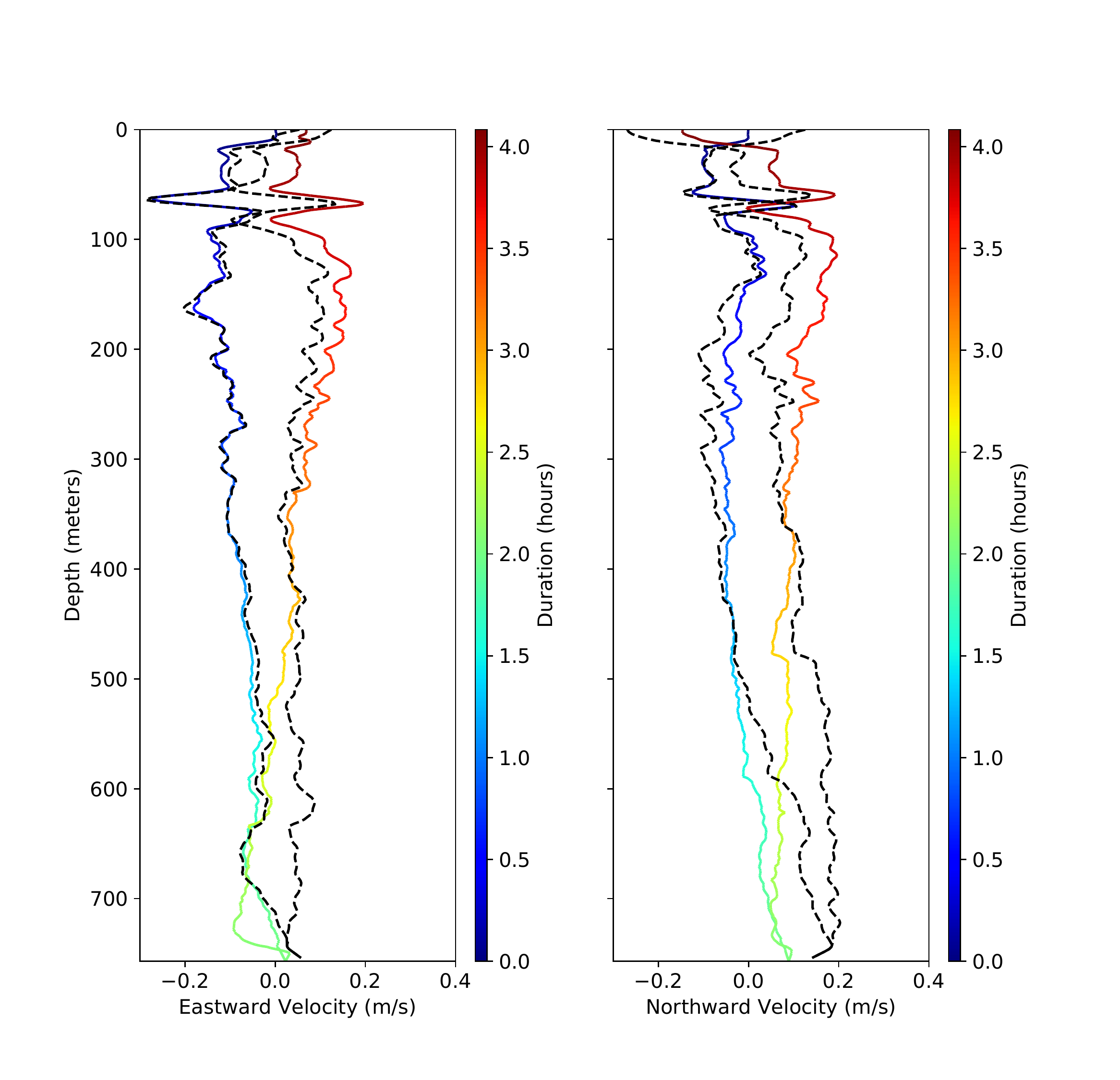}
  \caption{Comparison of the full results from the method of Section~\ref{sec:first} (black dashed) to that of Section~\ref{sec:nonlin} (colored).}
  \label{fig.comparison}
\end{figure}

\begin{figure}[]
  \includegraphics[width=\columnwidth]{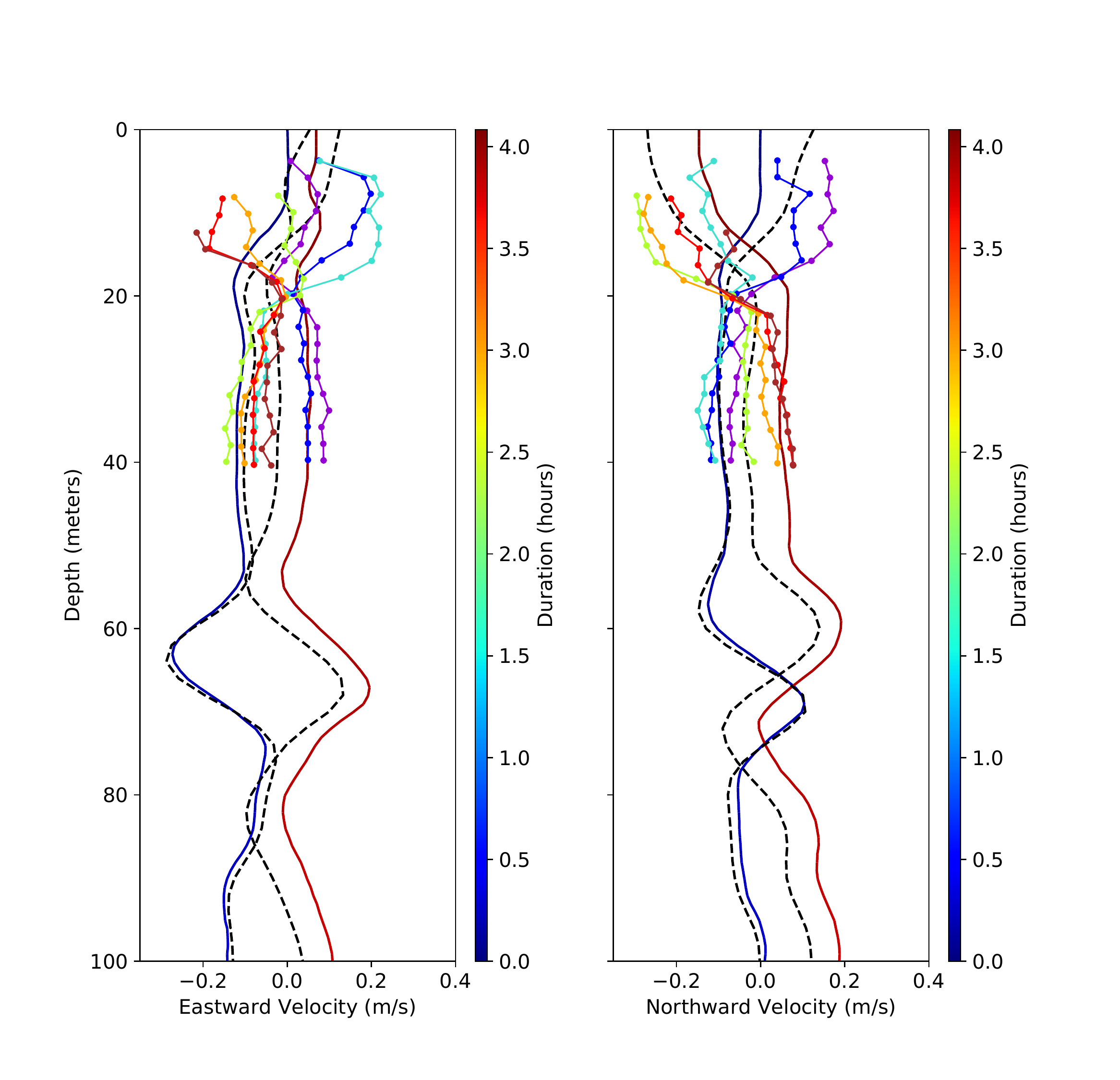}
  \caption{Zooming in on the upper 100m of the profiles in Figure \ref{fig.comparison}, we use current profile data from a 600 kHz ADCP on a mooring that was approx. 13 km away during this dive. Current data was collected hourly. The time evolution of the surface currents over the 5 hours during the dive, is captured by the color scale, starting with blue at the beginning of the dive and ending with red.}
  \label{fig.mooring}
\end{figure}

A comparison of the results of the two methods for the entire dive is shown in Figure \ref{fig.comparison}. Both approaches find highly correlated estimates of current profiles, which is reassuring, though there is a clear offset between the results that varies by depth. We believe that this is due in large part to the sensitivity of the state-space approach to the weighting applied to the different terms in the minimization, which are difficult challenging to vet without ground truth for the entire current profile.

The only ground truth available for this dive is for the upper 40 m, shown in Figure~\ref{fig.mooring}. These data are from an upward-facing 600 kHz moored ADCP that was approximately 13 km away during this dive. The time evolution of the surface currents, measured hourly, is captured by the same color scale used for the glider data, starting with blue at the beginning of the dive and ending with red.  Current data, as shown here, has not yet been corrected for the motion of the mooring. This correction and aggregate comparisons for other dives that are near moorings is the topic of future work.

\begin{figure}
  \includegraphics[width=\columnwidth]{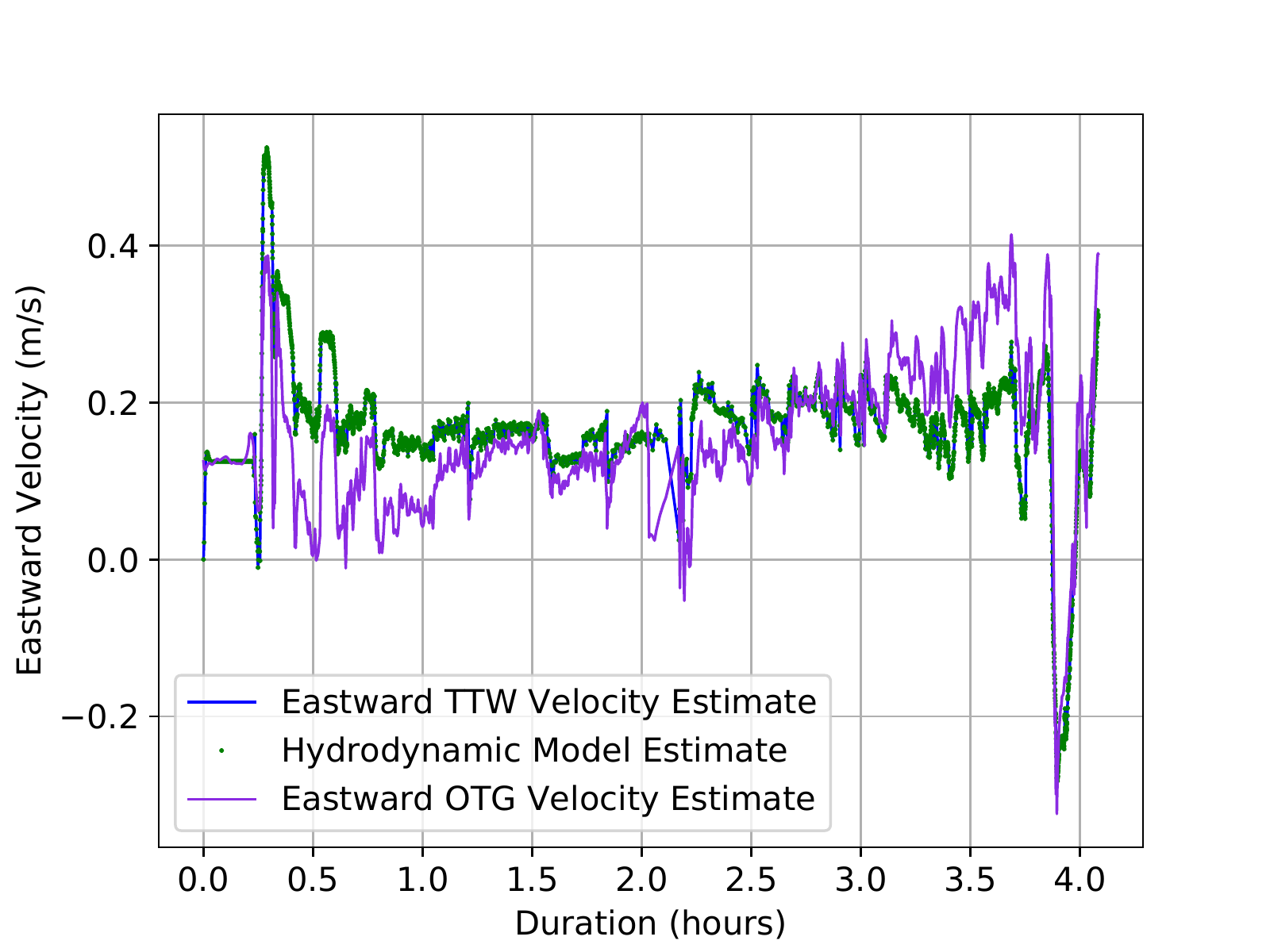}
  \includegraphics[width=\columnwidth]{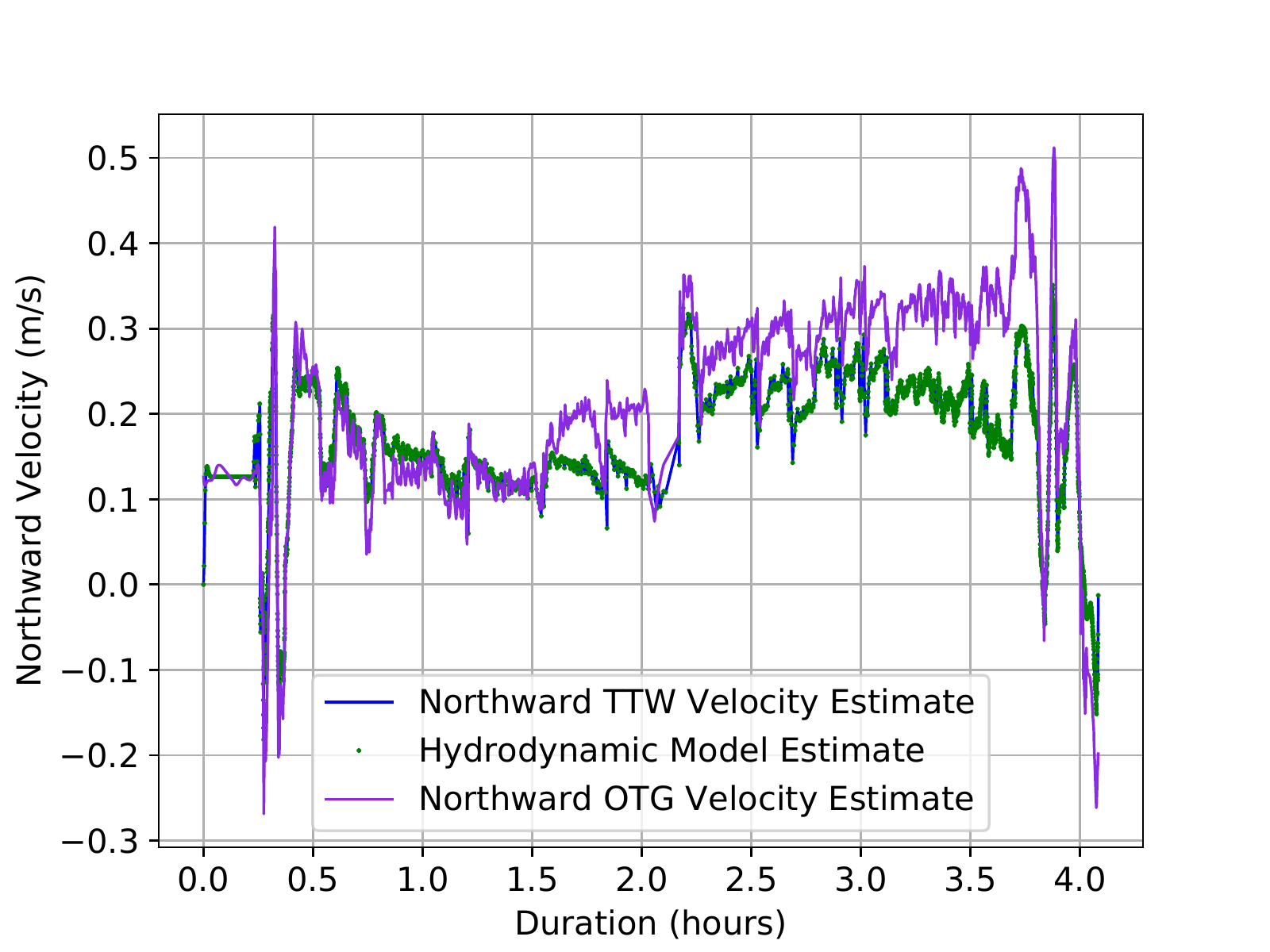}
  \caption{Updated velocities after processing the ADCP data. }
  \label{fig.updatedvel}
\end{figure}
\begin{figure}
    \centering
    \includegraphics[width=\columnwidth]{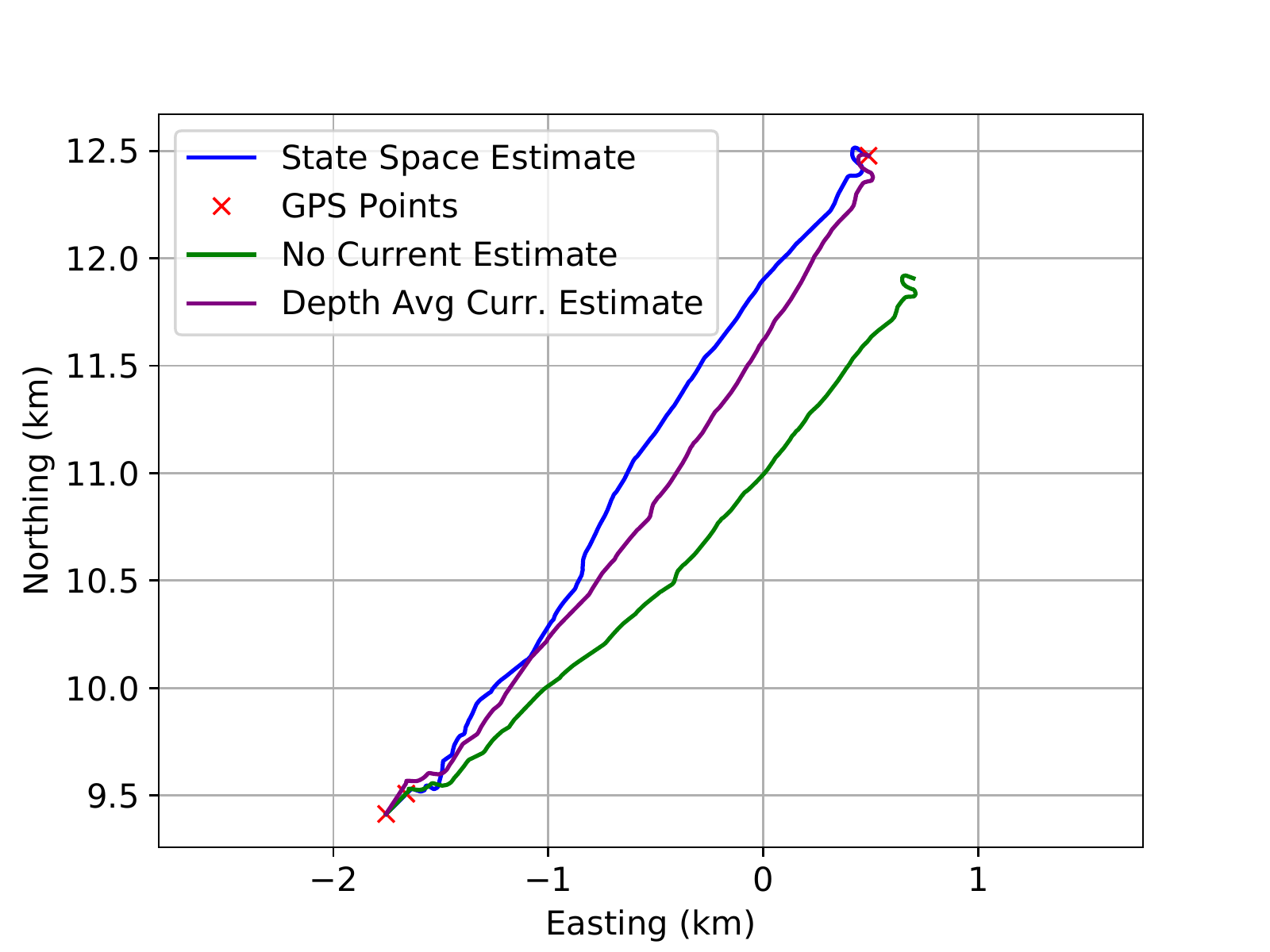}
    \caption{We compare position estimates from a glider trajectory computed without knowledge of current (dead reckoning), using the naive approach of applying a depth-averaged current uniformly across the dive, and using the ADCP-based current profile to inform the glider position throughout the dive. The ADCP-based correction not only shifts the trajectory, but also compresses the dive and stretches the climb, leading to a max horizontal offset of 449m.}
    \label{fig.xy}
\end{figure}

The second approach also produces updated velocity profiles, which are shown in Figure~\ref{fig.updatedvel} and the resulting glider trajectory, shown in Figure~\ref{fig.xy}. Comparing the x-y position estimates from a glider trajectory computed using the naive approach of uniformly applying a single depth-averaged current estimate across the dive versus using the ADCP-based current profile to inform the glider position throughout the dive, shows a dramatic difference. The ADCP-based correction not only shifts the trajectory, but also compresses the dive and stretches the climb. This leads to a maximum horizontal offset of 449m in this case. We look forward to using range measurements from the CANAPE moorings to both improve and validate these results in future work. 

\section{Conclusions and Future Work}

The main contribution of the work is to validate the idea of solving the state-space and current deconvolution problem using ADCP measurements. Two related approaches were developed, implemented, and compared; both made use of glider-estimated velocities, ADCP observations, and GPS fixes. The estimated current profiles were reasonably comparable to those obtained by a `ground truth' ADCP on a nearby mooring. 

The state-space approach in Section~\ref{sec:nonlin} is more informative than the first, since it also produces modified post-processed state estimates of the glider, see Figure~\ref{fig.updatedvel}. While there is no ground truth on the estimated velocities, the high degree of correlation between the deconvolution and inversion method for estimating current profiles validates the approach, and opens doors for future work using the state-space approach. 

First, the state-space approach builds an explicit connection between variables used during navigation to those informed by the ADCP. This is a promising development for ADCP-aided navigation. 

Second, the state-space framework allows a broad set of optimization-based tools to be applied to the deconvolution problem. In particular we can incorporate constraints~\cite{bell2009inequality}, robust losses~\cite{aravkin2014robust}, and nonsmooth regularization~\cite{aravkin2013sparse,aravkin2017generalized}, 
as well as singular state-space models that make use of these innovations~\cite{jonker2019fast}.
The state-space formulation also allows nonlinear models~\cite{bell2009inequality,aravkin2014optimization}, 
as needed by range measurements. All of these innovations make it possible to fuse information from noisy measurements along with statistical models that are robust to noisy data and constraints that incorporate prior knowledge. Efficient implementation of these formulations requires significant additional work in understanding the structure 
of the underlying linear algebra problem~\eqref{eq:KSjoint}.

\section*{Acknowledgements}

This work was completed as part of the Canada Basin Glider Experiment (CABAGE). Funding for this work was provided by the U.S. Office of Naval Research (ONR) through the Arctic and Global Predictions Program (Award \#N00014-16-2596, PI: Dr. Sarah Webster, APL-UW) and the Defense Research and Development Canada (DRDC) (Contract \#W7707-175902/001/HAL, PI: Dr. Sarah Webster, APL-UW).
Additional funding for CABAGE was provided by ONR Ocean Acoustics Program (Award \#N00014-17-1-2228, PI: Dr. Lora Van Uffelen, URI).

We would like to acknowledge Drs. Peter Worcester and Matthew Dzieciuch from Scripps Institution of Oceanography, who led the Canada Basin Acoustic Propagation Experiment (CANAPE). CANAPE provided all of the mooring infrastructure as well as all of the ice breaker logistics that enabled CABAGE.

This work could not have been completed without the engineers, technicians, and operators that made the Seaglider deployments possible. To that end, the authors would like to thank the Integrative Observational Platforms (IOP) Lab, in particular Craig Lee, who supplied the Seagliders used in the field experiment; Jason Gobat, the Senior Engineer; and Ben Jokinen, the Seaglider technician responsible for installing the ADCPs and preparing, testing, and launching the gliders. We are grateful for the support of Michael Flemming and Bill Kopplin, captain and current and former owners of the R/V Ukpik; the captain and crew of the USCGC Healy; and graduate student Wendy Snyder from University of Rhode Island, who assisted with the glider recovery. In addition, Jason Gobat and Geoff Shilling assisted with glider piloting.

\IEEEtriggeratref{8}

\bibliographystyle{IEEEtran}
\bibliography{sew,sew_addtl,references_sasha}

\end{document}

%% file: acronym_list.tex
%

\acrodef{2D}{two-dimensional}
\acrodef{3D}{three-dimensional}

\acrodef{ABE}{Autonomous Benthic Explorer}
\acrodef{ADCP}{acoustic Doppler current profiler}
\acrodef{ADV}{acoustic Doppler velocimeter}
\acrodef{Alvin}{Alvin}
\acrodef{AFRL}{Air Force Research Laboratory}
\acrodef{AHRS}{attitude and heading reference system}
\acrodef{Autosub}{Autosub}
\acrodef{AUG}{autonomous underwater glider}
\acrodef{AUV}{autonomous underwater vehicle}
\acrodef{AGU}{American Geophysical Union}
\acrodef{AON}{Arctic Observing Network}
\acrodef{AOPE}{Applied Ocean Physics and Engineering}
\acrodef{AS}{asymptotically stable}
\acrodef{ACFR}{Australian Centre for Field Robotics}
\acrodef{ASV}{autonomous surface vehicle}

\acrodef{BCC}{brightness constancy constraint}
\acrodef{BIO}{Bedford Institute of Oceanography}

\acrodef{CAD}{computer aided design}
\acrodef{CenSSIS}{Center for Subsurface Sensing and Imaging Systems}
\acrodef{CCD}{charge coupled device}
\acrodef{CI}{covariance intersection}
\acrodef{CG}{conjugate gradients}
\acrodef{COG}{course over ground}
\acrodef{CPU}{central processing unit}
\acrodef{CT}{continuous-time}
\acrodef{CB}{center of buoyancy}
\acrodef{CG}{center of gravity}
\acrodef{CGSN}{Canadian Gravity Standardization Net}
\acrodef{COM}{center of mass}
\acrodef{CO2}{carbon dioxide}
\acrodef{CPS}{control power supply}
\acrodef{CSAC}{chip-scale atomic clock}
\acrodef{CST}{{\it IEEE Transactions on Control Systems Technology}}
\acrodef{CSV}{Comma Separated Values}
\acrodef{CTD}{conductivity-temperature-depth}
\acrodef{CISE}{Computer and Information Science and Engineering}

\acrodef{DOF}{degree of freedom}
\acrodef{DoD}{Department of Defense}
\acrodef{DSL}{Deep Submergence Laboratory}
\acrodef{DT}{discrete-time}
\acrodef{DVL}{Doppler velocity log}
\acrodef{DPM}{digital panel meter}
\acrodef{DR}{dead reckoning}
\acrodef{DRI}{department research initiative}
\acrodef{DWH}{Deepwater Horizon}

\acrodef{ESDF}{exactly sparse delayed-state filter}
\acrodef{EM}{electromagnetic}
\acrodef{EKF}{extended Kalman filter}
\acrodef{EIF}{extended information filter}
\acrodef{ERC}{Engineering Research Center}
\acrodef{EPR}{East Pacific Rise}
\acrodef{EPSL}{{\it Earth and Planetary Science Letters}}

\acrodef{FastSLAM}{Factored Solution to SLAM}
\acrodef{FOG}{fiber optic gyro}
\acrodef{FOV}{field of view}
\acrodef{FFT}{fast Fourier transform}
\acrodef{FBD}{free body diagram}
\acrodef{FIRST}{For Inspiration and Recognition of Science and Technology}

\acrodef{GMRF}{Gaussian Markov random field}
\acrodef{GPS}{global positioning system}
\acrodef{GAS}{globally asymptotically stable}
\acrodef{GA}{geometric algebra}
\acrodef{GUI}{graphical user interface}

\acrodef{HUGIN}{HUGIN}
\acrodef{HMMV}{H\r{a}kon Mosby Mud Volcano}
\acrodef{HOV}{human occupied vehicle}
\acrodef{HROV}{hybrid remotely operated vehicle}
\acrodef{HTF}{Hydrodynamics Test Facility}

\acrodef{ICRA}{{\it IEEE International Conference on Robotics and Automation}}
\acrodef{IF}{information filter}
\acrodef{IFREMER}{French Institute for the Research and Exploitation of the Sea}
\acrodef{IFE}{Institute for Exploration}
\acrodef{INU}{inertial navigation unit}
\acrodef{INS}{inertial navigation system}
\acrodef{IR}{infrared}
\acrodef{IMU}{inertial measurement unit}
\acrodef{INS}{inertial navigation system}
\acrodef{IBCAO}{International Bathymetric Chart of the Arctic Ocean}
\acrodef{IROS}{{\it IEEE International Conference on Robotics and Intelligent Systems}}
\acrodef{iUSBL}{inverted USBL}
\acrodef{OWTTiUSBL}{One-Way Travel-Time inverted USBL}

\acrodef{Jason}{Jason}
\acrodef{JFR}{{\it Journal of Field Robotics}}
\acrodef{JdF}{Juan de Fuca}
\acrodef{JHU}{Johns Hopkins University}
\acrodef{JPEG}{Joint Photographic Experts Group}

\acrodef{KF}{Kalman filter}
\acrodef{KYP}{Kalman-Yakubovich-Popov}
\acrodef{KISS}{Keck Institute for Space Studies}

\acrodef{LADCP}{lowered acoustic Doppler current profiler}
\acrodef{LBL}{long-baseline}
\acrodef{LCM}{lightweight communication and marshaling}
\acrodef{LG}{linear Gaussian}
\acrodef{LKY}{Lefschetz-Kalman-Yakubovich}
\acrodef{LMedS}{least median of squares}
\acrodef{LLE}{Linear Lyapunov Equation}
\acrodef{LQR}{Linear Quadratic Regulation}
\acrodef{LQGR}{Linear Quadratic Gaussian Regulation}
\acrodef{LSSL}{Louis S. St Laurent}
\acrodef{LRAUV}{long-range AUV}

\acrodef{MAP}{maximum \emph{a posteriori}}
\acrodef{MBARI}{Monterey Bay Aquarium Research Institute}
\acrodef{MBN}{mosaic-based navigation}
\acrodef{MEF}{Main Endeavor Field}
\acrodef{MIT}{Massachusetts Institute of Technology}
\acrodef{MLE}{maximum likelihood estimate}
\acrodef{MRF}{Markov random field}
\acrodef{MIZ}{Marginal Ice Zone}
\acrodef{MATE}{Marine Advanced Technology Education}
\acrodef{MCR}{Mid-Cayman Rise}
\acrodef{MEMS}{micro-electrical-mechanical systems}
\acrodef{MBARI}{Monterey Bay Aquarium Research Institute}
\acrodef{MIMO}{multiple-input, multiple-output}
\acrodef{MOR}{mid-ocean ridge}
\acrodef{MVCO}{Martha's Vineyard Coastal Observatory}
\acrodef{MCL}{mission critical level}
\acrodef{MO}{mission objective}
\acrodef{MSD}{mass spring damper}
\acrodef{NavEst}{navigation estimation}
\acrodef{NDSEG}{National Defense Science and Engineering Graduate}
\acrodef{NEES}{normalized estimation error squared}
\acrodef{NDSF}{National Deep Submergence Facility}
\acrodef{NMEA}{National Marine Electronics Association}
\acrodef{NSF}{National Science Foundation}
\acrodef{NTP}{network time protocol}
\acrodef{NAS}{National Academy of Sciences}
\acrodef{NDSF}{National Deep Submergence Facility}
\acrodef{NHS}{Natick High School}
\acrodef{NLO}{nonlinear observer}
\acrodef{NSTA}{National Science Teachers Association}
\acrodef{NSIDC}{National Snow and Ice Data Center}
\acrodef{NUI}{Nereid under-ice}
\acrodef{OS}{operating system}
\acrodef{OWTT}{one-way travel time}
\acrodef{ODE}{ordinary differential equation}
\acrodef{OL}{Second-Order, Open-Loop Observer}
\acrodef{ONR}{Office of Naval Research}
\acrodef{OOI}{Ocean Observing Initiative}
\acrodef{OWTT}{one-way travel time}

\acrodef{PC}{personal computer}
\acrodef{PCB}{printed circuit board}
\acrodef{PDE}{partial differential equation}
\acrodef{PPS}{pulse per second}
\acrodef{PEG}{parameter error gain}
\acrodef{PHF}{Piccard Hydrothermal Field}
\acrodef{PNAS}{{\it Proceedings of the National Academy of Sciences}}
\acrodef{ppb}{part-per-billion}
\acrodef{ppm}{part-per-million}
\acrodef{PR}{Positive Real}
\acrodef{PROV}{polar remotely operated vehicle}
\acrodef{PI}{principal investigator}


\acrodef{RAM}{random access memory}
\acrodef{RANSAC}{random sample consensus}
\acrodef{RDI}{RD Instruments}
\acrodef{REMUS}{Remote Environmental Monitoring Unit}
\acrodef{RLG}{ring laser gyroscope}
\acrodef{RF}{radio frequency}
\acrodef{RMS}{Royal Mail Steamship}
\acrodef{ROM}{range only measurement}
\acrodef{ROV}{remotely operated vehicle}
\acrodef{RTC}{real-time clock}
\acrodef{ROV}{remotely operated vehicle}
\acrodef{ROVER}{remotely operated vehicle environmental research}
\acrodef{RNS}{Regional Scale Node}
\acrodef{RSS}{{\it Robotics: Science and Systems}}
\acrodef{RTS}{Rauch-Tung-Striebel}

\acrodef{SeaBED}{SeaBED}
\acrodef{SEIF}{sparse extended information filter}
\acrodef{SIFT}{scale invariant feature transform}
\acrodef{SLAM}{simultaneous localization and mapping}
\acrodef{SNAME}{The Society of Naval Architects and Marine Engineers}
\acrodef{SSD}{sum of squared differences}
\acrodef{SFM}{structure-from-motion}
\acrodef{SO}{special orthogonal}
\acrodef{SOA}{state of the art}
\acrodef{SPR}{Strictly Positive Real}
\acrodef{STEM}{Science, Technology, Engineering, and Mathematics }
\acrodef{SSF}{Summer Student Fellow}
\acrodef{SVO}{Scaler, Velocity Observer}
\acrodef{SVD}{singular value decomposition}
\acrodef{SVP}{sound velocity profile}
\acrodef{SPURS}{Salinity Processes in the Upper Ocean Regional Study}
\acrodef{TDMA}{time division multiple access}
\acrodef{TRL}{technology readiness level}
\acrodef{TRO}{{\it IEEE Transactions on Robotics}}
\acrodef{TJTF}{thin junction-tree filter}
\acrodef{TTL}{transistor-transistor logic}
\acrodef{TXCO}{temperature compensated crystal oscillator}
\acrodef{3D}{three dimensional}
\acrodef{TWTT}{two-way travel time}

\acrodef{USBL}{ultra-short-baseline}
\acrodef{UUV}{unmanned underwater vehicle}
\acrodef{UTC}{Coordinate Universal Time}
\acrodef{UDP}{User Datagram Protocol}
\acrodef{UKF}{Unscented Kalman Filter}
\acrodef{UV}{Underwater Vehicle}
\acrodef{UNOLS}{University National Oceanographic Laboratory System}
\acrodef{USGS}{United States Geological Survey}
\acrodef{USNA}{United States Naval Academy}
\acrodef{UNCLOS}{United Nations Convention on the Law of the Sea}
\acrodef{VAN}{visually augmented navigation}
\acrodef{VNL}{vision numerical library}
\acrodef{VTK}{The Visualization Toolkit}
\acrodef{VIGA}{vehicle induced gravimeter acceleration}

\acrodef{WHOI}{Woods Hole Oceanographic Institution}


\acrodef{YIP}{Young Investigator Program}


%% file: commands.tex
\newcommand{\preind}[3]{\;{{\tiny{#1}}\atop{\tiny{#2}}}\hspace{-0.05in}#3}
\newcommand{\laplace}{\mathcal{L}}
\newcommand{\rb}[1]{\raisebox{-1.5ex}[0pt]{#1}}

\newcommand{\sentry}[0]{{\it Sentry }}
\newcommand{\nereus}[0]{{\it Nereus }}
\newcommand{\jason}[0]{{\it Jason }}
\newcommand{\alvin}[0]{{\it Alvin }}
\newcommand{\iver}[0]{{\it Iver2 }}
\newcommand{\tioga}[0]{{\it Tioga }}

\newcommand{\order}[1]{\ensuremath{\mathcal{O}(10^{#1})}}

\newcommand{\deepGlider}{Deepglider}

\newcommand{\standardTraj}{conventional\xspace}
\newcommand{\USBLTraj}{deep-profiling\xspace}

\newcommand{\effFactor}{R}
\newcommand{\effFactorDive}{R_{dive}}

\newcommand{\TDive}{T_{dive}}
\newcommand{\Tconv}{T_{c}} 
\newcommand{\Tdeep}{T_{d}} 

\newcommand{\profileDepth}{z}
\newcommand{\profileHeight}{\Delta_z}
\newcommand{\diveSpeed}{U}
\newcommand{\vertSpeed}{W}
\newcommand{\diveAngle}{\theta}
\newcommand{\pumpEnergy}{E_{pump}}
\newcommand{\batteryEnergy}{B}
\newcommand{\pHotel}{P_{hotel}}
\newcommand{\pNav}{P_{nav}}
\newcommand{\diveEnergyConventional}{E_{dive,o}}
\newcommand{\diveEnergyUSBL}{E_{dive}}
\newcommand{\xenduranceConventional}{T_o}
\newcommand{\xenduranceUSBL}{T}


\newcommand{\vehicle}[0]{v}
\newcommand{\world}[0]{w}
\newcommand{\usbl}[0]{u}
\newcommand{\locallevel}[0]{n}  

\newcommand{\pos}[2]{\preind{#1}{}{\mathbf p}_{#2}}

\newcommand{\rot}[2]{\preind{#1}{#2}\mathbf{R}}

\newcommand{\upu}[0]{\pos{\usbl}{\mathrm{ASV}}}
\newcommand{\wpv}[0]{\pos{\locallevel}{\mathrm{ASV}}}
\newcommand{\vpw}[0]{\pos{\vehicle}{\world}}
\newcommand{\vuR}[0]{\rot{\vehicle}{\usbl}}
\newcommand{\uvR}[0]{\rot{\usbl}{\vehicle}}
\newcommand{\wvR}[0]{\rot{\locallevel}{\vehicle}}
\newcommand{\vwR}[0]{\rot{\vehicle}{\locallevel}}
\newcommand{\az}[0]{\alpha}
\newcommand{\el}[0]{\gamma}
\newcommand{\rng}[0]{\Gamma}

\newcommand{\x}[0]{x(t)}
\newcommand{\xdot}[0]{\dot{x}(t)}
\newcommand{\xddot}[0]{\ddot{x}(t)}
\newcommand{\xhat}[0]{\hat{x}(t)}
\newcommand{\xhatdot}[0]{\dot{\hat{x}}(t)}
\newcommand{\xhatddot}[0]{\ddot{\xhat}(t)}
\newcommand{\deltax}[0]{\Delta x(t)}
\newcommand{\deltaxdot}[0]{\Delta \dot{x}(t)}
\newcommand{\deltaxddot}[0]{\Delta \ddot{x}(t)}
\newcommand{\none}[0]{l_1 s_1 + l_2 s_2}
\newcommand{\ntwo}[0]{l_3 s_1 + l_4 s_2 + r}

\newcommand{\vel}[0]{v(t)}
\newcommand{\veldot}[0]{\dot{v}(t)}
\newcommand{\velhat}[0]{\hat{v}(t)}
\newcommand{\velhatdot}[0]{\dot{\hat{v}}(t)}
\newcommand{\veldelta}[0]{\Delta \vel}
\newcommand{\veldeltadot}[0]{\dot{\Delta \vel}}
\newcommand{\veldeltaq}[0]{\Delta v_q(t)}

\newcommand{\out}[0]{w(t)}
\newcommand{\outhat}[0]{\hat{w}(t)}
\newcommand{\outdelta}[0]{\Delta w(t)}

\newcommand{\invec}[0]{\bm{u}(t)}

\newcommand{\measnoise}[0]{n(t)}

\newcommand{\thrust}[0]{\tau(t)}
\newcommand{\mass}[0]{m}
\newcommand{\bouyancy}[0]{b}
\newcommand{\quaddrag}[0]{d_{Q}}
\newcommand{\lindrag}[0]{d_{L}}

\newcommand{\veli}[0]{v_i(t)}
\newcommand{\acceli}[0]{\dot{v}_i(t)}
\newcommand{\thrusti}[0]{\tau_i(t)}
\newcommand{\massi}[0]{m_i}
\newcommand{\alphai}[0]{\alpha_i}
\newcommand{\betai}[0]{\beta_i}
\newcommand{\mui}[0]{\mu_i}
\newcommand{\nui}[0]{\nu_i}
\newcommand{\bouyancyi}[0]{b_i}
\newcommand{\quaddragi}[0]{d_{Q_i}}
\newcommand{\lindragi}[0]{d_{L_i}}


\newcommand{\worldvel}[0]{ \preind{w}{}\dot{\bm{p}}}
\newcommand{\lblworldvel}[0]{ \preind{w}{}\dot{\bm{p}}_{l}}
\newcommand{\lblworldpos}[0]{ \preind{w}{}\bm{p}_{l}}

\newcommand{\beamvel}[0]{\bm{v}_{beam}(t)}
\newcommand{\dopinstvel}[0]{\preind{i}{}\dot{\bm{p}}_d(t)}
\newcommand{\dopworldvel}[0]{\preind{w}{}\dot{\bm{p}}_{d}(t)}
\newcommand{\dopworldposhat}[0]{\preind{w}{}\hat{\bm{p}}_d(t)}
\newcommand{\dopworldpos}[0]{\preind{w}{}\bm{p}_d(t)}
\newcommand{\dopworldposhatini}[0]{\preind{w}{}\hat{\bm{p}}_{d}(t_0)}
\newcommand{\dopvehvel}[0]{\preind{v}{}\dot{\bm{p}}_d(t)}

\newcommand{\cvec}[0]{\bm{c}}
\newcommand{\qvec}[0]{\bm{q}}
\newcommand{\qvecdot}[0]{\dot{\bm{q}}}

\newcommand{\xvec}[0]{\left[\begin{array}{{c}} \x \\ \xdot \end{array}\right]}
\newcommand{\xvecshort}[0]{{\bm{x}}(t)}
\newcommand{\xvecdot}[0]{\left[\begin{array}{{c}} \xdot \\ \xddot \end{array}\right]}
\newcommand{\xvecdotshort}[0]{\bm{\dot{x}}(t)}

\newcommand{\xhatvec}[0]{\left[\begin{array}{{c}} \xhat \\ \xhatdot \end{array}\right]}
\newcommand{\xhatvecshort}[0]{\bm{\hat{x}}(t)}
\newcommand{\xhatvecdot}[0]{\left[\begin{array}{{c}} \xhatdot \\ \xhatddot \end{array}\right]}
\newcommand{\xhatvecdotshort}[0]{\dot{\bm{\hat{x}}}(t)}

\newcommand{\deltaxvec}[0]{\left[\begin{array}{{c}} \Delta x(t) \\ \dot{\Delta x}(t) \end{array}\right]}
\newcommand{\deltaxvecshort}[0]{\Delta \bm{x}(t)}
\newcommand{\deltaxvecdot}[0]{\left[\begin{array}{{c}} \dot{\Delta x}(t) \\ \ddot{\Delta x}(t) \end{array}\right]}
\newcommand{\deltaxvecdotshort}[0]{\Delta \dot{\bm{x}}(t)}

\newcommand{\betavec}[0]{\left[\begin{array}{{c}} 0 \\ \beta \end{array}\right]}
\newcommand{\betavecshort}[0]{\bm{\beta}}

\newcommand{\alphavec}[0]{\left[\begin{array}{{c}} 0 \\ \alpha \end{array}\right]}
\newcommand{\alphavecshort}[0]{\bm{\alpha}}

\newcommand{\outvec}[0]{\bm{w}(t)}
\newcommand{\outvecshort}[0]{\bm{w}(t)} 
\newcommand{\outhatvecshort}[0]{\bm{\hat{w}}(t)}
\newcommand{\outtildevec}[0]{\bm{\tilde{w}}(t)}
\newcommand{\OUTtilde}[0]{\tilde{W}(t)}
\newcommand{\outdeltavecshort}[0]{\Delta \bm{w}(t)} 

\newcommand{\rvec}[0]{\left[\begin{array}{{c}} 0      \\ r_i(t)   \end{array} \right]}
\newcommand{\svec}[0]{\left[\begin{array}{{c}} s_1(t) \\ s_2(t) \end{array} \right]}
\newcommand{\nuvec}[0]{\left[\begin{array}{{c}} 0 \\ \nu \end{array}\right]}
\newcommand{\nuvecshort}[0]{\bm{\nu}}
\newcommand{\muvecshort}[0]{\bm{\mu}}
\newcommand{\muvec}[0]{\left[\begin{array}{{c}} 0 \\ \mu \end{array}\right]}

\newcommand{\measnoisevec}[0]{\bm{n}(t)}

\newcommand{\instframe}[0]{\preind{i}{}\bm{p}(t)}
\newcommand{\vehicleframe}[0]{\preind{v}{}\bm{p}(t)}
\newcommand{\localframe}[0]{\preind{l}{}\bm{p}(t)}
\newcommand{\worldframe}[0]{\preind{w}{}\bm{p}(t)}

\newcommand{\xq}[0]{\vel |\vel|}
\newcommand{\xqshort}[0]{x_q(t)}
\newcommand{\xhatq}[0]{\dot{\hat{x}}(t) |\dot{\hat{x}}(t)|}
\newcommand{\xhatqshort}[0]{\hat{x}_q(t)}
\newcommand{\deltaxq}[0]{\dot{\hat{x}}(t) |\dot{\hat{x}}(t)| - \dot{x}(t) |\dot{x}(t)|}
\newcommand{\deltaxqshort}[0]{\Delta \dot{x}_q(t)}

\newcommand{\Amatrix}[0]{\left[\begin{array}{{cc}} 0 & 1 \\ 0 & \mu \end{array}\right]}
\newcommand{\Ashort}[0]{\bm{A}}
\newcommand{\Cshort}[0]{\bm{C}}
\newcommand{\Pshort}[0]{\bm{P}}
\newcommand{\Lshort}[0]{\bm{L}}
\newcommand{\Qshort}[0]{\bm{Q}}
\newcommand{\eye}[0]{ I}
\newcommand{\momega}[0]{{\bm{m}}(i\omega)}
\newcommand{\ms}[0]{{\bm{m}}(s)}
\newcommand{\inmap}[0]{{\bm{b}}}
\newcommand{\outmap}[0]{\bm{c}}
\newcommand{\Rshort}[0]{\bm{R}}
\newcommand{\insttovehicle}[0]{\preind{v}{i}R}
\newcommand{\vehicletolocal}[0]{\preind{l}{v}R(t)}
\newcommand{\localtoworld}[0]{\preind{w}{l}R(t)}
\newcommand{\vehicletoworld}[0]{\preind{v}{w}R}

\newcommand{\cvecTtrans}[0]{\cvec\hspace{0.03in}^T}
\newcommand{\cvecT}[0]{\cvec}
\newcommand{\cvectrans}[0]{\cvec\hspace{0.03in}^T}
\newcommand{\bvecTtrans}[0]{\bvec\hspace{0.03in}^T}
\newcommand{\bvectrans}[0]{\bvec\hspace{0.03in}^T}
\newcommand{\bvecT}[0]{\bvec}
\newcommand{\qvectrans}[0]{\qvec\hspace{0.03in}^T}
\newcommand{\AT}[0]{A}
\newcommand{\QT}[0]{Q}
\newcommand{\PT}[0]{P}

\newcommand{\zetavec}[0]{\bm{\zeta}}

\newcommand{\gravposvec}[0]{\preind{w}{}\bm{p_g}}
\newcommand{\gravposx}[0]{\preind{w}{}p_{g_x}}
\newcommand{\gravposy}[0]{\preind{w}{}p_{g_y}}
\newcommand{\gravposz}[0]{\preind{w}{}p_{g_z}}
\newcommand{\gravposveclong}[0]{\left[ \begin{array}{c} \gravposx \\  \gravposy \\ \gravposz \end{array} \right]}
\newcommand{\gravposvecdot}[0]{\preind{w}{}\bm{\dot{p}_g}}
\newcommand{\gravposxdot}[0]{\preind{w}{}\dot{p}_{g_x}}
\newcommand{\gravposydot}[0]{\preind{w}{}\dot{p}_{g_y}}
\newcommand{\gravposzdot}[0]{\preind{w}{}\dot{p}_{g_z}}
\newcommand{\gravposvecddot}[0]{\preind{w}{}\bm{\ddot{p}_g}}
\newcommand{\gravposxddot}[0]{\preind{w}{}\ddot{p}_{g_x}}
\newcommand{\gravposyddot}[0]{\preind{w}{}\ddot{p}_{g_y}}
\newcommand{\gravposzddot}[0]{\preind{w}{}\ddot{p}_{g_z}}

\newcommand{\gravoffsetvec}[0]{\preind{g}{v}\bm{d}}
\newcommand{\gravoffsetx}[0]{\preind{g}{v}d_x}
\newcommand{\gravoffsety}[0]{\preind{g}{v}d_y}
\newcommand{\gravoffsetz}[0]{\preind{g}{v}d_z}
\newcommand{\gravoffsetveclong}[0]{\left[ \begin{array}{c} \gravoffsetx \\  \gravoffsety \\ \gravoffsetz \end{array} \right]}

\newcommand{\Rvehtograv}[0]{\preind{g}{v}R}
\newcommand{\hdg}[0]{\phi}
\newcommand{\pitch}[0]{\theta}
\newcommand{\roll}[0]{\psi}
\newcommand{\hdgdot}[0]{\dot{\phi}}
\newcommand{\pitchdot}[0]{\dot{\theta}}
\newcommand{\rolldot}[0]{\dot{\psi}}
\newcommand{\hdgddot}[0]{\ddot{\phi}}
\newcommand{\pitchddot}[0]{\ddot{\theta}}
\newcommand{\rollddot}[0]{\ddot{\psi}}

\newcommand{\vehicleposvec}[0]{\preind{w}{}\bm{p_v}}
\newcommand{\vehicleposx}[0]{\preind{w}{}p_{v_x}}
\newcommand{\vehicleposy}[0]{\preind{w}{}p_{v_y}}
\newcommand{\vehicleposz}[0]{\preind{w}{}p_{v_z}}
\newcommand{\vehicleposveclong}[0]{\left[ \begin{array}{c} \vehicleposx \\  \vehicleposy \\ \vehicleposz \end{array} \right]}
\newcommand{\vehicleposxdot}[0]{\preind{w}{}\dot{p}_{v_x}}
\newcommand{\vehicleposydot}[0]{\preind{w}{}\dot{p}_{v_y}}
\newcommand{\vehicleposzdot}[0]{\preind{w}{}\dot{p}_{v_z}}
\newcommand{\vehicleposxddot}[0]{\preind{w}{}\ddot{p}_{v_x}}
\newcommand{\vehicleposyddot}[0]{\preind{w}{}\ddot{p}_{v_y}}
\newcommand{\vehicleposzddot}[0]{\preind{w}{}\ddot{p}_{v_z}}

\newcommand{\vehiclevelvec}[0]{\left[ \begin{array}{c} \vehicleposxdot \\  \vehicleposydot \\  \vehicleposzdot \\ \hdgdot \\ \pitchdot \\ \rolldot \end{array} \right]}
\newcommand{\vehicleaccvec}[0]{\left[ \begin{array}{c} \vehicleposxddot \\  \vehicleposyddot \\  \vehicleposzddot \\ \hdgddot \\ \pitchddot \\ \rollddot \end{array} \right]}